\documentclass[12pt]{amsart}

\usepackage{amsmath}

\usepackage{amscd}

\usepackage{graphicx}

\usepackage{latexsym}

\usepackage{hyperref}

\usepackage{times}

\oddsidemargin .25in \theoremstyle{plain}

\newtheorem{Thm}{Theorem}

\newtheorem{Cor}[Thm]{Corollary}

\newtheorem{Prop}[Thm]{Proposition}

\newtheorem{Def}[Thm]{Definition}

\newtheorem{Rem}[Thm]{Remark}

\newcommand{\bfc}{{\mathbb{C}}}

\newcommand{\OB}{\mathfrak{{ob}}}

\def\v{\vskip.12in}

\begin{document}

\title[]{Explicit horizontal open books on some plumbings}

\author{Tolga Etg\"u}

\author{Burak Ozbagci}

\address{Department of Mathematics \\ Ko\c{c} University \\ Istanbul, Turkey}

\email{tetgu@ku.edu.tr and bozbagci@ku.edu.tr}

\address{ Current: School of Mathematics, Georgia Institute
of Technology, Atlanta,
Georgia, USA}

\email{bozbagci@math.gatech.edu}

\subjclass[2000]{57R17, 53D10, 32S55, 32S25 }

%57R57, 53D35, 57R95}

\date{\today}

\thanks{T.E. was partially supported by the Turkish Academy of Sciences;
B.O. was partially supported by the Turkish Academy of Sciences and by
the NSF Focused Research Grant FRG-024466.}

\begin{abstract}

We describe explicit open books on arbitrary plumbings of oriented
circle bundles over closed oriented surfaces. We show that, for a
non-positive plumbing, the open book we construct is horizontal and
the corresponding compatible contact structure is also horizontal
and Stein fillable. In particular, we describe horizontal open books
on some Seifert fibered 3--manifolds. As another application we
describe horizontal open books isomorphic to Milnor open books for
some complex surface singularities. Moreover we give examples of
tight contact 3--manifolds supported by planar open books. As a
consequence the Weinstein conjecture holds for these tight contact
structures \cite{ach}.

\end{abstract}

\maketitle %\setcounter{section}{-1}

\section{Introduction}

In this article using 3--dimensional simple surgery techniques we
first construct explicit open books on oriented circle bundles over
closed oriented surfaces. Then we construct open books on plumbings
of circle bundles according to a graph by appropriately gluing the
open books we constructed for the circle bundles involved in the
plumbing. (We will use the word graph for a connected graph with at
least two vertices.) An open book on a circle bundle or more
generally on a plumbing of circle bundles is called
\emph{horizontal} if its binding is a collection of some fibers and
its pages are transverse to the fibers. Here we require that the
orientation induced on the binding by the pages coincides with the
orientation of the fibers induced by the fibration.

We will call a plumbing graph \emph{non-positive} if the sum of the
degree of the vertex and the Euler number of the bundle
corresponding to that vertex is non-positive for every vertex of the
graph. We prove that the open book we construct on a circle bundle
with negative Euler number or on a non-positive plumbing of circle
bundles is horizontal. It turns out that the contact structure
compatible with this open book is also \emph{horizontal}, i.e. the
contact planes (possibly after an isotopy of the contact structure)
are transverse to the fibers. Furthermore we show that the monodromy
of this horizontal open book is given by a product of right-handed
Dehn twists along disjoint curves. Consequently, by a theorem of
Giroux \cite{gi}, our horizontal open book is compatible with a
Stein fillable contact structure. Recall that Eliashberg and Gromov
\cite{eg} proved that every fillable contact structure is tight. As
a first application of our construction we describe explicit
horizontal open books on some Seifert fibered 3--manifolds.

A contact 3--manifold $(Y, \xi)$ is said to be \emph{Milnor
fillable} if it is contactomorphic to the contact boundary of an
isolated complex surface singularity $(X, x)$. The germ  $(X, x)$ is
called a \emph{Milnor filling} of $(Y, \xi)$. An analytic function
$f : (X, x) \to (\mathbb{C}, 0)$ with an isolated singularity at $x$
defines an open book decomposition of $Y$ which is called a
\emph{Milnor open book}. It was shown in \cite{cnp} that any Milnor
open book of a Milnor fillable contact 3--manifold $(Y, \xi)$ is
horizontal and compatible with the contact structure $\xi$. We would
like to point out that a Milnor fillable contact structure is Stein
fillable and hence tight.

A Milnor fillable contact 3--manifold is obtained as a plumbing of
some circle bundles over surfaces which is topologically completely
determined by the minimal good resolution of the singularity. As a
second application we describe explicit open books isomorphic to
Milnor open books for some complex surface singularities.

On the other hand, in \cite{ga}, using 4--dimensional symplectic
handle attachments, Gay gives a construction of open books on
plumbings of circle bundles. It turns out that the open book we
construct on a given plumbing is isomorphic to Gay's open book,
showing in particular that the open book in his construction can be
made horizontal for non-positive plumbings. If a plumbing graph is
not non-positive, then there are binding components in our open book
which are oriented opposite to the fiber orientation, and hence the
open book fails to be horizontal.

Finally we give examples of planar open books whose compatible
contact structures are Stein fillable and hence tight. As a
consequence the Weinstein conjecture holds for these tight contact
structures, since the conjecture is proved for every contact
structure compatible with a planar open book (cf. \cite{ach}).
Planar open books compatible with some Stein fillable contact
structures were also constructed independently by Sch\"{o}nenberger
in his thesis \cite{sc}, using completely different techniques.

 \v \noindent {\bf {Acknowledgement}}: We are grateful to
John Etnyre and Andr\'{a}s Stipsicz  for their critical comments on
the first draft of this paper that helped us to improve our earlier
results. We would also like to thank Mohan Bhupal, Ferit \"Ozt\"urk,
and Meral Tosun for helpful conversations.

\section{Open book decompositions and contact structures}

We will assume throughout this paper that a contact structure
$\xi=\ker \alpha$ is coorientable (i.e., $\alpha$ is a global
1--form) and positive (i.e., $\alpha \wedge d\alpha >0 $ ). In the
following we describe the compatibility of an open book
decomposition with a given contact structure on a 3--manifold.

Suppose that for an oriented link $L$ in a 3--manifold $Y$ the
complement $Y\setminus L$ fibers over the circle as $\pi \colon Y
\setminus L \to S^1$ such that $\pi^{-1}(\theta) = \Sigma_\theta $
is the interior of a compact surface bounding $L$, for all $\theta
\in S^1$. Then $(L, \pi)$ is called an \emph{open book
decomposition} (or just an \emph{open book}) of $Y$. For each
$\theta \in S^1$, the surface $\Sigma_\theta$ is called a
\emph{page}, while $L$ the \emph{binding} of the open book. The
monodromy of the fibration $\pi$ is defined as the diffeomorphism of
a fixed page which is given by the first return map of a flow that
is transverse to the pages and meridional near the binding. The
isotopy class of this diffeomorphism is independent of the chosen
flow and we will refer to that as the \emph{monodromy} of the open
book decomposition.

An open book $(L, \pi)$ on a 3--manifold $Y$ is said to be
\emph{isomorphic} to an open book $(L^\prime, \pi^\prime)$ on a
3--manifold $Y^\prime$, if there is a diffeomorphism $f: (Y,L)  \to
(Y^\prime, L^\prime)$ such that $\pi^\prime \circ f = \pi$ on $Y
\setminus L$. In other words, an isomorphism of open books takes
binding to binding and pages to pages.

{\Thm [Alexander \cite{al}] Every closed and oriented 3--manifold
admits an open book decomposition.}

\vskip1ex

In fact, every closed oriented 3--manifold admits
a planar open book, which means that a page is a disk $D^2$
with holes \cite{rol}.
On the other hand, every closed
oriented 3--manifold admits a contact structure \cite{marti}.
So it seems natural to strengthen the contact condition
$\alpha \wedge d \alpha >0$ in the presence of an open book decomposition
on $Y$
by requiring that $\alpha>0 $ on the binding and $d\alpha>0 $ on
the pages.

{\Def \label{compatible} An open book decomposition of a 3--manifold
$Y$ and a contact structure $\xi$ on $Y$ are called
\emph{compatible} if $\xi$ can be represented by a contact form
$\alpha$ such that the binding is a transverse link, $d \alpha$ is a
symplectic form on every page and the orientation of the transverse
binding induced by $\alpha$ agrees with the boundary orientation of
the pages.}

\begin{Thm} [Giroux \cite{gi}] \label{giroux} Every contact 3--manifold
admits a compatible open book. Moreover two contact structures
carried by the same open book are isotopic.
\end{Thm}

In particular, in order to show that two contact structures on a
given 3--manifold are contactomorphic, it suffices to show that they
are carried by isomorphic open books. We refer the reader to
\cite{et2} and \cite{ozst} for more on the correspondence between
open books and contact structures.

\section{Explicit construction of open books} \label{construction}

We will assume throughout the paper that the circle bundles we
consider are oriented and the base space is a closed oriented
surface. An open book on a single circle bundle or on a plumbing of
circle bundles (according to a graph) is called \emph{horizontal} if
its binding is a collection of some fibers in the circle bundles,
its pages are transverse to the fibers of the circle bundles and the
orientation induced on the binding by the pages coincides with the
orientation of the fibers induced by the fibration.

We start with constructing
an open book on a circle bundle with negative Euler number.
We show that our construction yields a horizontal open book,
whose compatible contact structure is horizontal as well.
Then we generalize the construction to circle bundles
with non-negative Euler numbers, but the open
books we construct on those bundles are not horizontal.

\subsection{Construction of an open book on an oriented
circle bundle over an oriented surface}

\begin{figure}[ht]

  \begin{center}

     \includegraphics{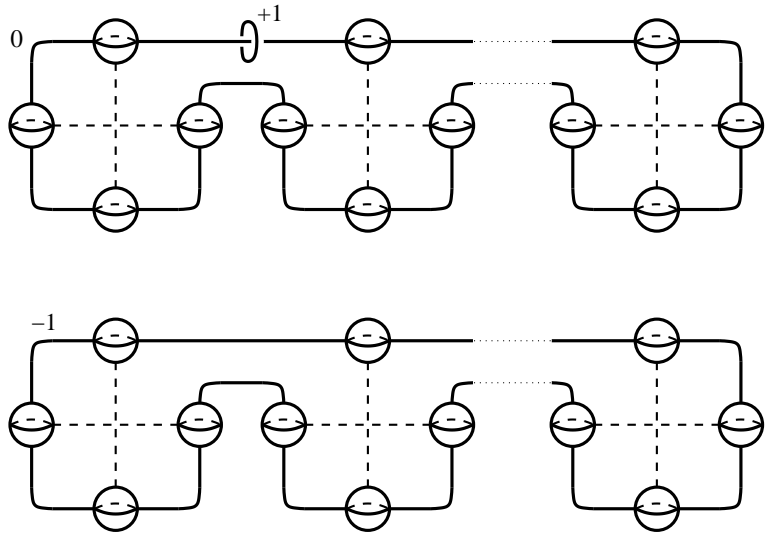}

   \caption{The diagram on top shows $D^2 \times\Sigma$ (whose boundary is $S^1 \times \Sigma$) along with a
$+1$--surgery curve transverse to $\Sigma$. One can visualize the surface $\Sigma$ as the disk
(bounded by the $0$--framed knot) union $2$-dimensional $1$--handles going over the $4$--dimensional $1$--handles.
We blow down the $+1$
curve to get a disk bundle over $\Sigma$ with Euler number $-1$. The
 boundary of this disk bundle at the bottom is a circle bundle over
 $\Sigma$ with Euler number $-1$. } \label{circbundle}

    \end{center}
  \end{figure}

Let $\Sigma$ denote a closed oriented surface and let $M$ denote the
trivial circle bundle over $\Sigma$, i.e., $M = S^1 \times \Sigma$.
We pick a circle fiber $K$ of $M$ and perform $+1$--surgery along
$K$. The resulting 3--manifold $M^\prime$ is a circle bundle over
$\Sigma$ whose Euler number is equal to $-1$. This can be verified
using Kirby calculus (cf. \cite{gs}), for example, as shown in
Figure~\ref{circbundle}.

In order to describe an open book on $M^\prime$ we first consider
the affect of removing a solid torus neighborhood $N= K \times
D^2$ of $K$ from $M$. By removing $N$ from $M$ we puncture once
each $\Sigma$ in $M = S^1 \times \Sigma$ to get $S^1 \times
\widetilde{\Sigma}$, where $\widetilde{\Sigma}= \Sigma \setminus
D^2$. Now we will glue a solid torus back to $S^1 \times
\widetilde{\Sigma}$ along its boundary torus $S^1 \times \partial
\widetilde{\Sigma}$ in order to perform our surgery.

\begin{figure}[ht]

  \begin{center}

     \includegraphics{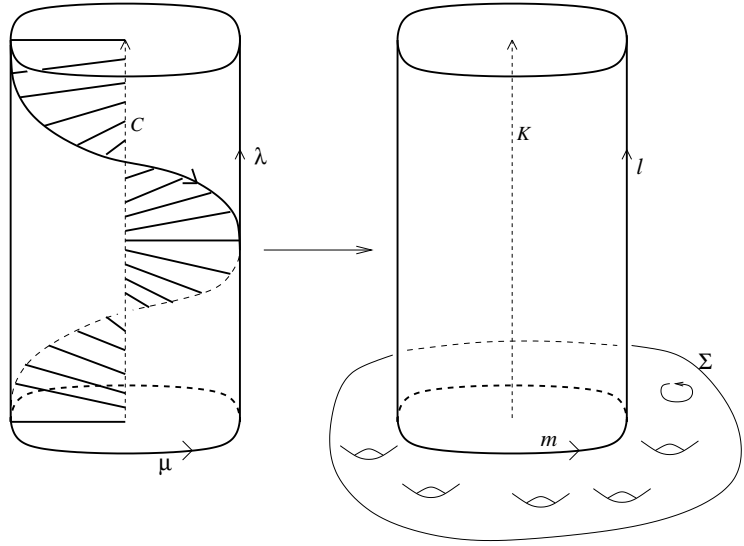}

   \caption{The surgery torus is obtained by
 identifying the top and the bottom of the cylinder shown on
the left-hand side.
On the right hand side we depict the complement obtained by
removing a neighborhood of a circle fiber from
$S^1 \times \Sigma$. } \label{surgerytorus}

    \end{center}

  \end{figure}

Consider the solid torus $S^1 \times D^2$ shown on the left-hand
side in Figure~\ref{surgerytorus}. Let $\mu$ and $\lambda$ be the
meridian and the longitude pair of $S^1 \times \partial D^2$. Let
$m$ and $l$ denote the meridian and the longitude pair in the
boundary of $S^1 \times \widetilde{\Sigma}$. Note that the base
surface is oriented (as we depicted in Figure~\ref{surgerytorus})
and the orientation induced on $ \partial  \widetilde{\Sigma}$ is
the opposite of the orientation of $m$.  We glue the solid torus
$S^1 \times D^2$ to $S^1 \times \widetilde{\Sigma}$ by an
\textit{orientation preserving} diffeomorphism $ S^1 \times
\partial D^2 \to S^1 \times \partial \widetilde{\Sigma}$
which sends $\mu$ to $m+l$ and $\lambda$ to $ l$. The resulting
3--manifold $M^\prime$ will be oriented extending the orientation on
$S^1 \times \widetilde{\Sigma}$ induced from $M$. We use orientation
preserving map rather than reversing because of our choice of
orientations above.

In Figure~\ref{surgerytorus}, we also depict a leaf (an annulus) of
a foliation on the solid torus $S^1 \times D^2$ that we will glue to
perform $+1$--surgery. Note that $S^1 \times D^2$ is a trivial
circle bundle over $D^2$ where the circle fibers are transverse to
the annuli foliation. The boundary of a leaf consists of the core
circle $C$ of $S^1 \times D^2$ and a $(1,-1)$--curve on $S^1 \times
\partial D^2$, i.e., a curve homologous to $ \mu - \lambda$. Each
leaf is oriented so that the induced orientation on the boundary of
a leaf is given as indicated in Figure~\ref{surgerytorus}. The
gluing diffeomorphism maps $\mu - \lambda$ onto $m$ so that by
performing the $+1$--surgery we also glue each annulus in the
foliation to a $\widetilde{\Sigma}$ in $S^1 \times
\widetilde{\Sigma}$ identifying the outer boundary component (i.e.,
the $(\mu - \lambda)$--curve) of the annulus with $\partial
\widetilde{\Sigma}$. Hence this construction yields an open book on
$M^\prime$ whose binding is $C$ (the core circle of the surgery
torus) and pages are obtained by gluing an annulus to each
$\widetilde{\Sigma}$ along $\partial \widetilde{\Sigma}$. Notice
that the pages will be oriented extending the orientation on
$\widetilde{\Sigma}$ induced from $\Sigma$. Finally we want to point
out that the core circle $C$ becomes an oriented fiber of the circle
fibration of $M^\prime$ over $\Sigma$.

Next we will describe the monodromy of this open book. In order to
measure the monodromy of an open book we should choose a flow which
is transverse to the pages and meridional near the binding. Take a
\emph{vertical} vector field pointing along the fiber direction in
$S^1 \times \widetilde{\Sigma}$ and extend it inside the surgery
torus (along every ray towards the core circle) by rotating
clockwise so that it becomes horizontal near the core circle.
Observe that since the vector field is horizontal near the binding,
the first return map of the flow generated by this vector field will
fix the points near the binding on any leaf. Now take a horizontal
arc (on a leaf) connecting the core circle to the other boundary of
that leaf. Then one can see that the flow will move the points of
this arc further to  the right if we move towards the boundary. Note
that the first return map will fix the points on
$\widetilde{\Sigma}$ as well as the points near the binding.
Combining this discussion we conclude that the first return map is
given by a right-handed Dehn twist along the core circle of the
leaf, which is indeed a curve parallel to the boundary of the open
book.

\begin{figure}[ht]

  \begin{center}

     \includegraphics{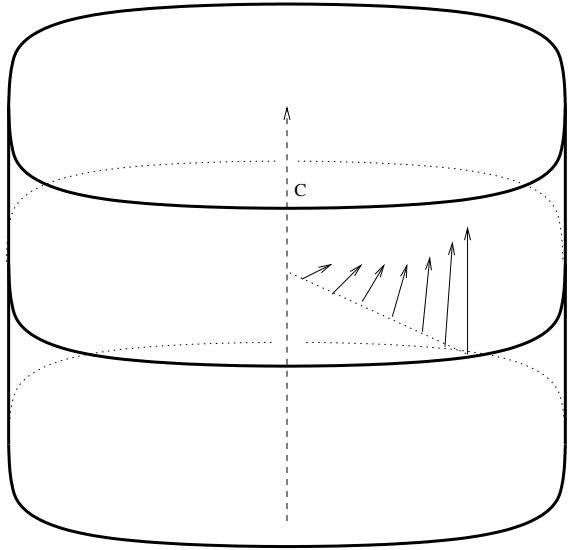}

   \caption{The vector field depicted along a ray towards
the core circle inside the surgery torus extending the
vertical one on $S^1 \times \widetilde{\Sigma}$.
The vector field becomes horizontal near the core (binding) along any ray.} \label{flow}

    \end{center}

  \end{figure}

In summary we constructed a horizontal open book on a circle bundle
over a genus $g$ surface $\Sigma$ with Euler number $-1$ by
performing a $+1$--surgery on a circle fiber of a trivial circle
bundle. It is clear that if we take $s$ disjoint fibers in $S^1
\times \Sigma$ and perform $+1$--surgery along each of these fibers,
we will get a horizontal open book on a circle bundle over $\Sigma$
with Euler number $e=-s < 0$. The binding of the resulting open book
will be the union of $s$ circle fibers and a page will be a genus
$g$ surface with $s$ boundary components. The monodromy will be the
product of right-handed Dehn twists along boundary components. Consequently
the compatible contact structure is Stein fillable (cf. \cite{gi}).
Thus we get

{\Prop \label{hor} There exists an explicit \emph{horizontal}
open book on a circle bundle with negative Euler number
which is compatible with a Stein fillable contact structure.}

\vspace{1ex}

Next we consider the contact structure compatible with our open book.

{\Prop \label{transverse} The contact structure compatible with the
open book we constructed on a circle bundle with negative Euler
number is \emph{horizontal}, i.e., the contact planes (possibly
after an isotopy of the contact structure) are positively transverse
to the fibers.}

\begin{proof}

By Lemma 3.5 in \cite{et2} we may assume that the contact planes are
arbitrarily close to the tangents of the pages away from the
binding. So clearly the contact planes are positively transverse to
the fibers away from the binding since the pages of our open book
are already positively transverse to the fibers.

Note that near a component $C$ of the binding we have explicitly
constructed the pages, and the fibers of the circle fibration can be
viewed as straight vertical lines in the solid cylinder on the
left-hand side in Figure~\ref{surgerytorus} before we identify the
top and the bottom. On the other hand, by the compatibility
condition, there are coordinates $(z, (r, \theta)) $ near $C$, where
$C$ is $\{r=0\}$ and a page is given by setting $\theta$ equal to a
constant such that the contact structure is given by the kernel of
the form $dz + r^2 d \theta$.
\begin{figure}[ht]

    \begin{center}

     \includegraphics{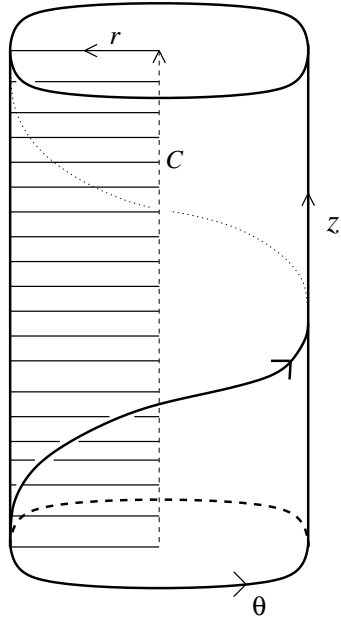}

   \caption{A circle fiber is given by a $(1,1)$--curve in a
neighborhood of a component $C$ of the binding.}  \label{twist}

    \end{center}
  \end{figure}
In these new coordinates the neighborhood of $C$ in
Figure~\ref{surgerytorus} is seen as in Figure~\ref{twist}, where
the annulus is straightened out and the circle fibers are wrapped
once around the cylinder in a right-handed manner. More precisely,
the tangent vector to an (oriented) fiber is given by
$\partial_{\theta} +
\partial_z$.
%Observe that the contact planes will be horizontal at the center $C$
%of the solid cylinder (cf. Figure~\ref{twist}) and they will rotate
%\emph{clockwise} along any horizontal ray starting at the center.
%Note that in these coordinates
%$$ \{ \partial_r, -r^2 \partial_z +
%\partial_\theta  \} $$ is an oriented
%basis of the contact plane $\ker (dz + r^2 d \theta )$.
Hence one can see that the contact planes  will remain positively
transverse to the circle fibers in a neighborhood of $C$, since
$$ (dz + r^2 d\theta) (\partial_{\theta} + \partial_z) = r^2 +1 > 0 .$$
\end{proof}

{\Rem In fact, the same argument can be generalized to prove that
the contact structure compatible with any horizontal open book on a
circle bundle is horizontal. Any horizontal contact structure on a
circle bundle is universally tight by Lemma 3.9 in \cite{ho}}

\vskip1ex

Note that the same construction of an open book will work when we
perform $-1$--surgeries along fibers, but in that case the
orientation on the binding induced by the pages will be the opposite
of the fiber orientation. Thus the open book we construct using
$-1$--surgeries will not be horizontal. Moreover we will get
left-handed Dehn twists along boundary components instead of
right-handed Dehn twists. Nevertheless, by taking  $s$ disjoint
fibers in $S^1 \times \Sigma$ and performing a $-1$--surgery along
each of these fibers, we will get an open book on a circle bundle
over $\Sigma$ with Euler number $e= s > 0$.

Finally, to obtain an open book on a circle bundle with zero Euler
number we first observe that by
performing a $+1$--surgery on a fiber and a $-1$--surgery in another
fiber does not change the topology of the 3--manifold which fibers
over a surface. Hence by performing a canceling pair of a $\pm 1$
surgeries on distinct fibers, we can construct an open book on a
circle bundle with zero Euler number
with two binding components, where the monodromy is given by
a right-handed Dehn twists along one boundary component
and a left-handed Dehn twist along the other.

\subsection{Construction of an open book on a plumbing of
oriented circle bundles over oriented surfaces}

We will construct open books for plumbings of circle bundles. We start
with describing an open book on the plumbing of two circle bundles. We obtain
our open book by suitably gluing the open books we constructed on the
circle bundles in the previous section.

Let $M_i$ denote a circle bundle over a closed surface $\Sigma_i$
with Euler number $e_i$ for $i = 1,2$. In the previous section we
constructed an open book on $M_i$. We observed that by
performing a $+1$--surgery on a fiber and a $-1$--surgery in another
fiber does not change the topology of the 3--manifold which fibers
over a surface. Thus we can assume that $M_i$ is obtained from $S^1
\times \Sigma$ by $\pm1$--surgeries along (disjoint fibers) such
that there are as many fibers as we wish on which we perform
$+1$--surgeries.

\begin{figure}[ht]

  \begin{center}

     \includegraphics{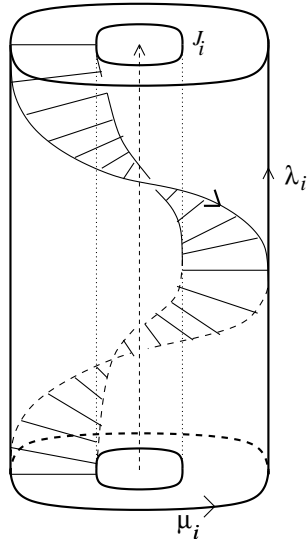}

   \caption{Removing a smaller neighborhood of the binding} \label{halfannulus}

    \end{center}

  \end{figure}

Fix one of the circles where we performed a $+1$--surgery to obtain
$M_i$. Consider the surgery process we described in the previous
section. Take a smaller solid torus neighborhood $J_i$ (see
Figure~\ref{halfannulus}) of the surgery torus we glued in to
perform $+1$--surgery. To obtain a plumbing of $M_1$ and $M_2$ we
remove $J_i$ from $M_i$ and identify boundary tori $\partial J_1$,
$\partial J_2$ by an orientation reversing diffeomorphism $
\partial J_1 \to \partial J_2$ taking $\mu_1$ to $ \lambda_2$
and $ \lambda_1$ to $
\mu_2$. Here as in the previous section $(\mu_i, \lambda_i)$
denote the meridian and the longitude pair for $J_i$. Note that
this diffeomorphism will take $\mu_1- \lambda_1$ to $ -(\mu_2 -
\lambda_2)$.

The point here is that when we remove $J_i$ from $M_i$ we remove a
(small) annulus neighborhood of a binding component from every page.
The remaining part of the annulus in the pages will intersect the
boundary torus $\partial J_i$ as a curve homologous to $\mu_i
-\lambda_i$. Hence when we plumb the circle bundles we actually glue
the pages (after removing a binding component) of the open books on
$M_i$ to obtain an open book on the plumbing. The new page will be
obtained from the pages of the open books on $M_1$ and $M_2$ by
gluing together one boundary component of the open book on $M_1$
with one boundary component of the open book on $M_2$ (both of which
were induced by $+1$--surgery). Note that the glued up page will be
oriented extending the orientations of the pages that we glue
together. Moreover the pages will be transverse
to the fibers by construction.

To calculate the monodromy of the open book on the plumbed manifold
we just observe that the flows we used to calculate monodromies of
the open books on each piece in the plumbing will also glue together:
We can assume that the flow has a constant slope $1$ (given by tangents
to the curves $\mu_i + \lambda_i$) on each boundary torus so
that the gluing map will take one to the other to
define a flow on the closed
3-manifold obtained by plumbing. When we
glue the two flows, since each flow will have the affect of moving the points
on the half annuli (see Figure~\ref{halfannulus}) to the right by $180$-degrees
we get a right-handed
Dehn twist along the annulus obtained by gluing the two half annuli.

Next we describe an explicit open book on arbitrary plumbings of
circle bundles according to a weighted graph. Suppose that we are
given a plumbing according to a weighted graph with $k$ vertices
such that each vertex is decorated with a pair of integers $(e_i,
g_i)$, where $e_i$ is the Euler number and $g_i$ is the genus of the
base surface $\Sigma_i$ of the $i^{th}$ circle bundle $M_i$ in the
plumbing. Let $d_i$ denote the degree of the $i^{th}$ vertex.

First assume that $e_i \geq 0$. Then we can view $M_i$ as obtained
from $S^1 \times \Sigma_i$ by performing $-1$--surgeries on $e_i$
distinct fibers. To be able to plumb $M_i$ to other circle bundles
$d_i$ times, we need to use $d_i$ fibers with $+1$--surgeries. So we
apply extra $+1$--surgeries on $d_i$ fibers as well as
$-1$--surgeries on some other $d_i$ fibers not to change the
topology of $M_i$. Since we will use all the $+1$--surgeries in the
plumbing process we will end up with $e_i + d_i$ fibers on which we
performed $-1$--surgeries.

Now assume that $e_i < 0 $. Then we can view $M_i$ as obtained from
$S^1 \times \Sigma_i$ by performing $+1$--surgeries on $-e_i$
distinct fibers. If we further assume that $d_i > -e_i$  we perform
$d_i+e_i$ extra $+1$ and $-1$--surgeries and use all the
$+1$--surgeries for the plumbing to end up with $d_i+e_i$ fibers on
which we applied $-1$--surgeries. The only remaining case is $d_i
\leq -e_i$ when $e_i <0$. In this case we can use $d_i$ of these
$+1$--surgeries in the plumbing and we will end up with $-e_i-d_i$
fibers on which we applied $+1$--surgeries.

Combining all the discussion above, a page of the open book on this
plumbing will be obtained by gluing together surfaces $F_i$ of genus
$g_i$ with $|e_i + d_i|$ boundary components according to the given
graph. (Note, however, that we will not use these boundary
components in the gluing.) Each edge in the graph will become a
1--handle (a neck) in the page connecting the surfaces corresponding
to the vertices connected by that edge. The monodromy of the open
book will be given as a product of one boundary-parallel Dehn twist
for each boundary component of the page and one right-handed Dehn
twist along the cocore circle of each 1--handle coming from an edge.
If $e_i+ d_i > 0 $ ($e_i+ d_i < 0$, resp.) for some $i$, then the
boundary parallel Dehn twists in $F_i$ will be left-handed
(right-handed, resp.). Here we will assume that \emph{there is at
least one vertex in a plumbing graph such that $e_i+ d_i$ is
non-zero}, to avoid the case of empty binding which we do not want
to allow in the definition of an open book. In particular, note that
if the plumbing graph is non-positive then our open book will be
horizontal since we do not use any $-1$--surgeries in that case.
Moreover the compatible contact structure is Stein fillable (cf.
\cite{gi}) since the monodromy of the open book for a non-positive
plumbing is a product of right-handed Dehn twists. Hence we get

{\Thm \label{horiz} There exists an explicit \emph{horizontal} open book on a
\emph{non-positive} plumbing of circle bundles over surfaces which
is compatible with a Stein fillable contact structure.}

\vspace{1ex}

Furthermore  we have

{\Prop \label{tans} The contact structure compatible with the open
book we constructed on a \emph{non-positive} plumbing of circle
bundles over surfaces is \emph{horizontal}, i.e., the contact planes
(possibly after an isotopy of the contact structure) are positively
transverse to the fibers.}

\begin{proof} In Proposition~\ref{transverse}, we proved that
the contact structure compatible with the
open book we constructed on a circle bundle with negative
Euler number is horizontal. Recall that the open book on a plumbing of two circle bundles
is obtained by removing a neighborhood of a binding component
from each open book and gluing along the resulting boundary tori.
The result follows, by the same argument used in Proposition~\ref{transverse},
since we proved that the resulting open book on the plumbing
is horizontal.
\end{proof}

\vspace{1ex}

{\Rem \label{horizon} In fact, the same argument can be generalized
to prove that the contact structure compatible with any horizontal
open book on a plumbing is horizontal.}

\vskip1ex

\begin{figure}[ht]

  \begin{center}

     \includegraphics{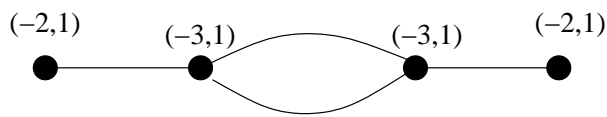}

   \caption{The integer weights at each vertex are the Euler number and the
genus of the base of the corresponding circle bundle, respectively.} \label{graph_ex1}

    \end{center}

  \end{figure}

\noindent {\bf{Example 1. }}Consider the 3--manifold $Y_1$ obtained
by plumbing circle bundles over tori according to the graph in
Figure~\ref{graph_ex1}. Degree-one vertices have $e_i + d_i = -1 $ but the other two
vertices have degree equal to three and $e_j+d_j = 0$.

\begin{figure}[ht]

  \begin{center}

     \includegraphics{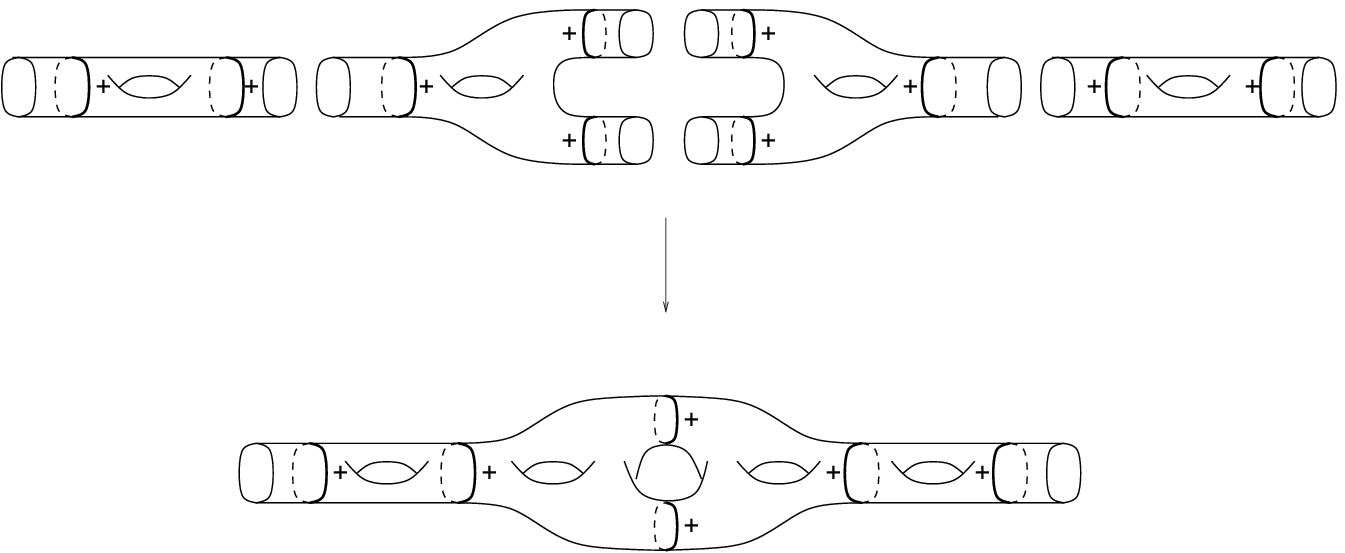}

   \caption{The open book $\OB_1$} \label{ex1_open_book}

    \end{center}

  \end{figure}

Therefore the
horizontal open book $\OB_1$ we construct according to the recipe
above will have a twice punctured (one puncture for each of the
degree-one vertices in the graph) genus 5 surface as a page and the
product of 2 right-handed Dehn twists along 2 disjoint separating
curves, 2 right-handed Dehn twists along 2 disjoint non-separating
curves and one right-handed Dehn twist parallel to each boundary
component as its monodromy (see Figure~\ref{ex1_open_book}). The
open book $\OB_1$ is an explicit horizontal open book on $Y_1$.

\vspace{1ex}

{\Rem In \cite{ga}, using 4--dimensional symplectic handle
attachments, Gay also gives a construction of open books on
plumbings of circle bundles. It turns out that the page and the
monodromy of the open book we construct on a given plumbing
coincides with the one given in \cite{ga}. This shows that the open
book we construct on a plumbing is (abstractly) isomorphic to Gay's
open book, showing that the open book in his construction can be
made horizontal for a non-positive plumbing.}

\vspace{1ex}

{\Cor A Seifert fibered 3--manifold $Y$ with an orientable base of
genus $g$ and Seifert invariants  (see \cite{ozst} for conventions) $$\{g,n;r_1, r_2, \cdots, r_k \}
$$ admits a horizontal open book if
$n +k \leq 0$. If $g=0$, this open book can be chosen to be planar.}

\begin{proof} By \cite{or}, any Seifert fibered 3--manifold $Y$ is isomorphic
(as a 3--manifold with $S^1$--action) to the boundary of a
4--manifold with $S^1$--action obtained by equivariant plumbing of
disk bundles according to a weighted star-shaped graph as in
Figure~\ref{seifert}, where
$$[b_{i1}, b_{i2}, \dots , b_{is_i}]=b_{i1}-
\frac{1}{b_{i2}-\frac{1}{\ddots - \frac{1}{b_{is_i}}}}$$ is the
unique continued fraction expansion of
$1/ r_i$ with each $b_{ij} \geq 2$.
Since each $b_{ij} \geq 2$, every vertex except for the central one
satisfies the non-positivity assumption. If the central vertex also
satisfies $n + k \leq 0$, then the Seifert fibered 3--manifold
admits a horizontal open book.

\begin{figure}[ht]

  \begin{center}

     \includegraphics{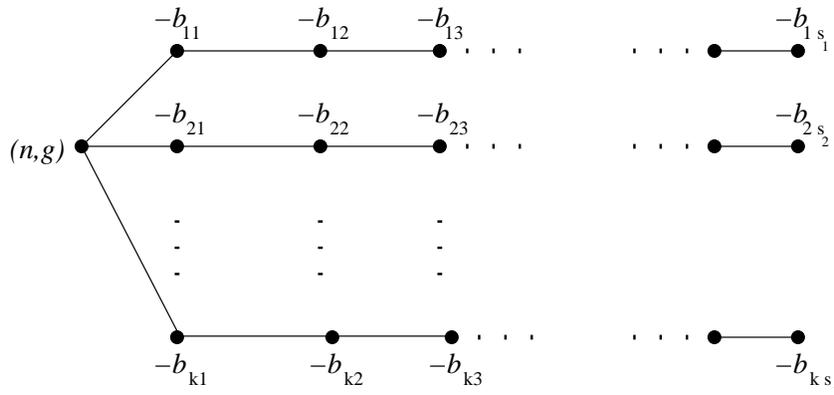}

   \caption{The plumbing graph of the Seifert fibered 3--manifold $Y$ with Seifert
invariants $\{g,n;r_1, r_2, \cdots, r_k \} .$ All bundles except the one at the central vertex are over spheres.}

\label{seifert}

    \end{center}

  \end{figure}

\end{proof}

\section{Milnor open books} \label{milopen}

\noindent

Let $(X,x)$ be an isolated complex-analytic singularity. Given a
local embedding of $(X,x)$ in $(\bfc^N, 0)$. Then a small sphere
$S^{2N-1}_{\epsilon} \subset \bfc^N$ centered at the origin intersects $X$ transversely,
and the complex hyperplane distribution $\xi$ on $Y=X\cap
S^{2N-1}_{\epsilon}$ induced by the complex structure on $X$ is a contact structure.
For sufficiently small
radius $\epsilon$ the contact manifold $(Y, \xi)$ is independent of
the embedding and $\epsilon$ up to contactomorphism and this
contactomorphism type is called the \emph{contact boundary} of
$(X,x)$. A contact manifold $(Y^\prime, \xi^\prime)$ is said to be
\emph{Milnor fillable} and the germ $(X, x)$ is called a
\emph{Milnor filling} of $(Y^\prime, \xi^\prime)$ if $(Y^\prime,
\xi^\prime)$ is contactomorphic to the contact boundary $(Y, \xi)$
of $(X, x)$. Even though these definitions are valid in all
dimensions, we will focus on surface singularities and their
boundaries of real dimension 3. Note that any surface singularity
can be normalized without changing the contact boundary and normal
surface singularities are isolated.

It is a well-known result of Grauert \cite{gr} that an oriented
3--manifold has a Milnor fillable contact structure if and only if
it can be obtained by plumbing oriented circle bundles over surfaces
according to a weighted graph with negative definite intersection
matrix. On the other hand, a recent discovery (cf. \cite{cnp}) is
that any 3--manifold admits at most one  Milnor fillable contact
structure, up to contactomorphism. This result is obtained by using
Milnor open books.

{\Def Given an analytic function $f: (X,x) \rightarrow (\bfc , 0)$
with an isolated singularity at $x$, the open book decomposition
$\OB_f$ of the boundary $Y$ of $(X,x)$ with binding $L=Y \cap
f^{-1} (0)$ and projection $ \pi=\frac{f}{|f|} : Y\setminus L \to
S^1 \subset \bfc$ is called the \emph{Milnor open book} induced by
$f$. }

Naturally, one can talk about Milnor open books on any Milnor
fillable manifold. Milnor open books have certain features that are
used in the proof of the uniqueness result mentioned above: they are
all compatible with the natural contact structure on the boundary
and they are horizontal when considered on the plumbing description
of the Milnor fillable 3--manifold.

{\Rem By Remark \ref{horizon} any Milnor fillable contact structure
is horizontal.}

\vskip1ex

On the other hand, if $Y$ is obtained by plumbing $k$ circle
bundles, we obtain a $k$--tuple ${\bf{n}}=(n_1, \dots , n_k)$ for
each horizontal open book of $Y$, where $n_i$ is the number of
distinct fibers in the $i^{th}$ circle bundle which appear in the
binding. Each $n_i$ is nonnegative by definition. According to
Proposition 4.6 in \cite{cnp}, if each $n_i$ is positive and the
plumbing graph has a nondegenerate intersection matrix, then
$\bf{n}$ uniquely determines the horizontal open book up to
isomorphism. This proposition can be applied to Milnor fillable
3--manifolds since they are obtained by plumbing circle bundles
according to negative definite intersection matrices. Also note
that, by Giroux correspondence, isomorphic open books are compatible
with contactomorphic contact structures.

As a result, given a 3--manifold $Y$ with a Milnor fillable contact
structure $\xi$, to show that any other Milnor fillable contact
structure $\xi'$ on $Y$ is contactomorphic to $\xi$, it is suffices
to show the existence of a $k$--tuple $\bf{n}$ of positive integers
such that for any Milnor filling $(X,x)$ of $(Y, \xi')$ there is a
holomorphic function $f: (X,x) \rightarrow (\bfc , 0)$ with an
isolated singularity at $x$ and a Milnor open book $\OB_f$ which
generates the $k$--tuple $\bf{n}$. It turns out that such a
$k$--tuple $\bf{n}$ exists, and in fact, every $\bf{n}$ that
satisfies the following two conditions works (see Theorem 4.1 in
\cite{cnp}):

\begin{enumerate}

\item $n_i \geq d_i + 2g_i$ for each $i$, where $g_i$ is the genus
of the base of the $i^{th}$ circle bundle and $d_i$ is the degree of
the $i^{th}$ vertex, \item there exists a $k$--tuple ${\bf{m}}$ of
nonnegative integers such that ${\bf{I \cdot m = - n}} \ ,$ where
${\bf{I}}$ is the (negative definite) intersection matrix of the
weighted plumbing graph, and ${\bf{m}}$ and ${\bf{n}}$ are
considered as column vectors.

\end{enumerate}

\section{Explicit Milnor open books}\label{ex}

Let $(Y, \xi)$ be a Milnor fillable contact 3--manifold. Then $Y$
can be obtained by plumbing circle bundles. Let $G$ be the weighted
graph (with $k$ vertices and a negative definite intersection matrix
$\bf{I}$) of such a plumbing, and $d_i$ and $(e_i, g_i)$ denote the
degree and the weight of the $i^{th}$ vertex of $G$ respectively.

{\Prop  \label{case} Suppose that $$e_i + 2d_i + 2g_i \leq 0$$ holds
for every vertex of $G$. Then the open book of $Y$ we construct in
Section \ref{construction} is horizontal and isomorphic to a Milnor
open book, for any Milnor filling of $Y$. Moreover the monodromy of
this horizontal open book is given by a product of
\emph{right-handed} Dehn twists along disjoint curves.}

\begin{proof}

This open book is horizontal and its monodromy is given by a product
of right-handed Dehn twists along disjoint curves simply because the
plumbing is non-positive. The entries of the $k$--tuple $\bf{n}$ it
generates are $n_i = |e_i + d_i|$. The inequality $e_i + 2d_i + 2g_i
\leq 0$ implies that $n_i \geq d_i +2g_i$, and the system ${\bf{I
\cdot m = - n}}$ has a solution, namely ${\bf{m}}=(1, \dots , 1) $.
Therefore the conditions (1) and (2) in Section~\ref{milopen} are
satisfied for the vector ${\bf n}$ implying the isomorphism between
this horizontal open book and a Milnor open book.

\end{proof}

\vspace{1ex}

\begin{figure}[ht]

  \begin{center}

     \includegraphics{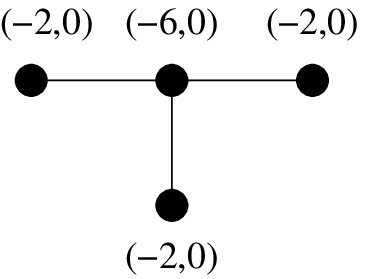}

   \caption{} \label{d4}

    \end{center}

  \end{figure}

\noindent{\bf{Example 2.}} Consider the 3--manifold $Y_2$ obtained
by plumbing circle bundles over spheres according to the graph in
Figure~\ref{d4}.

\begin{figure}[ht]

  \begin{center}

     \includegraphics{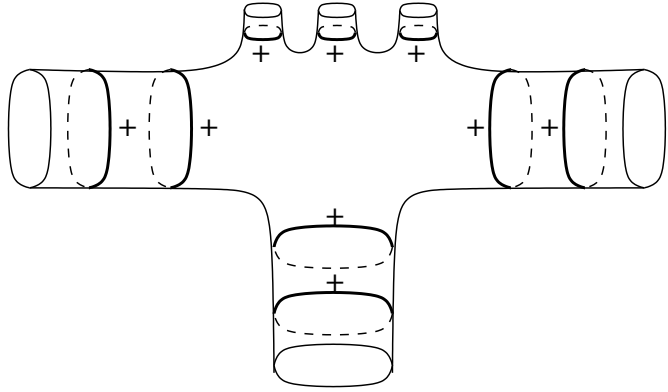}

   \caption{The open book $\OB_2$} \label{d4_open_book}

    \end{center}

  \end{figure}

The vertices with degree one have $e_i + d_i = -1 $ and the fourth
vertex has $e_4+d_4 = -3$. Note that the intersection matrix
$${{\bf I}} = \left(
  \begin{array}{cccc}
    -2 & 0 & 0& 1 \\
    0 & -2 & 0 & 1 \\
    0 & 0 & -2 & 1 \\
    1 & 1 & 1 & -6 \\
  \end{array}
\right)
$$
of the plumbing graph is negative definite and hence $Y_2$ admits a
Milnor fillable contact structure. The horizontal open book $\OB_2$
we construct for this non-positive plumbing according to the recipe
in Section~\ref{construction} will have a 6--times punctured sphere
as a page. The monodromy of $\OB_2$ is depicted in
Figure~\ref{d4_open_book}.

The quadruple generated by this horizontal open book is
${\bf{n}}=(1,1,1,3)$.  The system ${\bf{I \cdot m = - n}}$ has a
solution ${\bf{m}}=(1,1,1,1)$ and as a consequence the horizontal
open book $\OB_2$ is isomorphic to a Milnor open book for any Milnor
filling of $Y_2$. In particular, $\OB_2$ is compatible with the
unique Milnor fillable contact structure on $Y_2$.
Note that $Y_2$ is the Seifert fibered
manifold with invariants $\{0,-6; 1/2, 1/2, 1/2 \}$.

\section{Tight planar open books}\label{tightplanar}

In \cite{et1}, Etnyre showed that any overtwisted contact structure
on a closed 3--manifold is compatible with a planar open book. He
also provided the first obstructions for fillable contact structures to
be compatible with planar open books.

Recall that for any given contact $1$--form $\alpha$ there is a
unique vector field $X_\alpha$
defined by the conditions:
$$\alpha (X_\alpha) =1 \; \mbox{and}\; \iota_{X_\alpha} d \alpha = 0 .$$
The vector field $X_\alpha$ is called the Reeb vector field of
$\alpha$. Weinstein conjecture says that any Reeb vector field has a
closed orbit. For example, every component of the binding of an open
book is a closed orbit of the Reeb vector field of any contact
$1$--form $\alpha$ with $\alpha> 0$ on the binding and $d\alpha>0 $
on the pages of the open book. Recently, the conjecture is proved
for every contact structure compatible with a planar open book (cf.
\cite{ach}).

In this section we give examples of horizontal, Stein fillable
contact structures compatible with planar open books. As a
consequence the Weinstein conjecture holds for these tight contact
structures.

{\Cor Let $Y$ be a 3--manifold obtained by a plumbing of circle
bundles over spheres according to a tree. Then the open book we
construct on this plumbing is planar. If we assume that the plumbing
is non-positive then our planar open book is horizontal and
compatible with a Stein fillable and horizontal contact structure on
$Y$. Suppose furthermore that the inequality
$$e_i +2d_i \leq 0$$ holds for every vertex of the graph. Then our
planar horizontal open book is compatible with the
unique Milnor fillable contact structure on $Y$.}

\begin{proof}

Consider the open book of $Y$ constructed in
Section~\ref{construction}. The first claim immediately follows by the
fact that the plumbing is
according to a tree and the circle bundles involved are over
spheres. The second claim follows from the first and Theorem~\ref{horiz}.
To prove the last claim we just observe that
the condition $e_i +2d_i \leq 0$ trivially implies that the
plumbing is negative, i.e. $e_i+d_i < 0$. Hence the plumbing is
negative-definite, Y is Milnor fillable and also the open book is
horizontal. The open book is isomorphic to a Milnor open book for
any Milnor filling of $Y$ by Proposition~\ref{case}.

\end{proof}

\begin{figure}[ht]

  \begin{center}

     \includegraphics{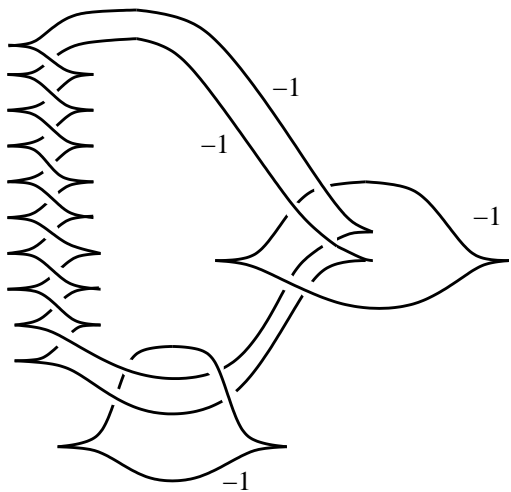}

   \caption{A surgery diagram for the unique Milnor fillable contact structure
on $Y_2$. (The framings are relative to the contact framing.)}

\label{milfill}

    \end{center}

  \end{figure}

\begin{figure}[ht]

  \begin{center}

     \includegraphics{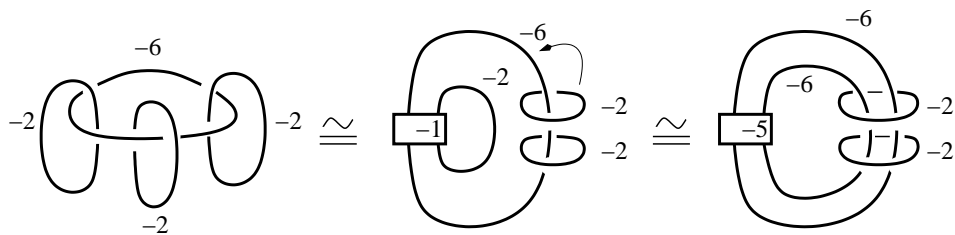}

   \caption{Sliding a 2-handle.} \label{262}

    \end{center}

  \end{figure}

For example, the open book $\OB_2$ depicted in
Figure~\ref{d4_open_book} is a planar horizontal open book
compatible with the unique Milnor fillable contact structure on
$Y_2$. In fact, a surgery diagram of this contact structure is given
in Figure~\ref{milfill}. We obtained this diagram by comparing our
construction with the one described in \cite{sc}. In
Figure~\ref{262}, we illustrate how to slide a 2-handle to show that
$Y_2$ is diffeomorphic to the 3--manifold in Figure~\ref{milfill}.

\section{An overtwisted example}\label{ot}

\vspace{1ex}

\begin{figure}[ht]

  \begin{center}

     \includegraphics{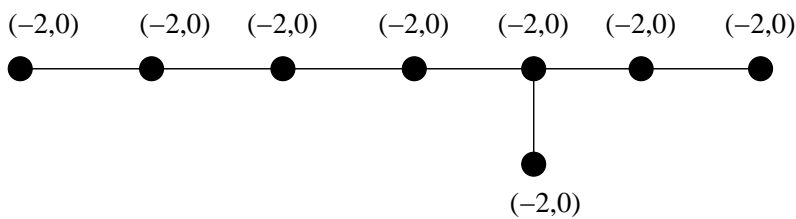}

   \caption{Negative definite $E_8$--plumbing} \label{e8}

    \end{center}

  \end{figure}

\begin{figure}[ht]

  \begin{center}

     \includegraphics{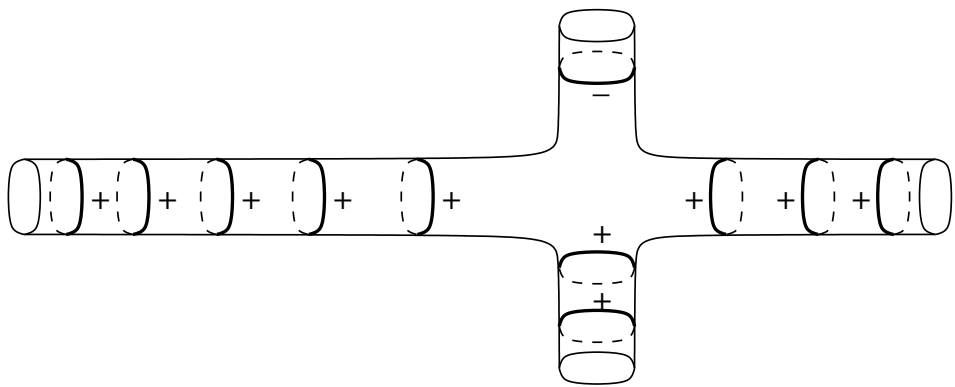}

   \caption{ $\OB_3$ : An overtwisted open book on the Poincar\'{e} homology sphere.}

\label{e8_open_book}

    \end{center}

  \end{figure}

\noindent{\bf Example 3.} Consider the singularity given by the
equation
$$ x^2 +y^3+z^5 = 0 $$
in $\bfc^3$. From the minimal resolution we obtain the \emph{dual
graph} depicted in Figure~\ref{e8} and hence the 3--manifold $Y_3$
obtained by plumbing circle bundles of Euler number $-2$ over
spheres according to this graph is the boundary of the singularity.
Note that the  intersection matrix of
this weighted graph is negative definite and the open book
$\OB_3$ we construct on $Y_3$ (which is just is the Poincar\'{e} homology sphere $\Sigma(2,3,5)$)
is depicted in Figure~\ref{e8_open_book}.

\begin{figure}[ht]

  \begin{center}

     \includegraphics{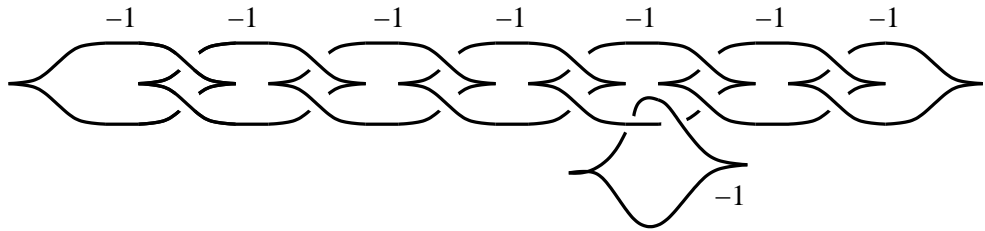}

   \caption{A surgery  diagram of the unique Milnor fillable contact structure on the Poincar\'{e} homology sphere. (The framings are relative to the contact framing.)}

\label{sigma235}

    \end{center}

  \end{figure}

We note that the contact structure compatible with $\OB_3$ is
overtwisted since
the Poincar\'{e} homology sphere does not admit any planar open book
compatible with its unique tight contact structure \cite{et1}. On the other hand,
a contact surgery diagram of the unique tight contact structure
on the Poincar\'{e} homology sphere is
depicted in Figure~\ref{sigma235}.
This contact structure is clearly the unique Milnor fillable
contact structure on  $\Sigma(2,3,5)$.

\end{document}